\documentclass[10pt]{amsart}
\usepackage{amsfonts,amssymb,amscd,amsmath,enumerate,verbatim,calc}
\usepackage{color,soul}
\usepackage{mathtools}
\usepackage{graphicx}
\everymath{\displaystyle}
\usepackage{fancyhdr}
\usepackage{tikz}
\usepackage{url}
\usetikzlibrary{shapes.geometric, arrows}
\tikzstyle{arrow} = [thick,->,>=stealth]
\tikzstyle{process} = [rectangle, minimum width=4cm, minimum height=2cm, text centered, text width=2.5cm, draw=green, fill=green!10]
\theoremstyle{plain}

\theoremstyle{definition}

\begin{document}
\title{H*-Normal Spaces and Some Related Functions}
\author{Neeraj kumar Tomar}
\email{neer8393@gmail.com}
\address{Department of Applied Mathematics, Gautam Buddha University, Greater Noida, Uttar Pradesh 201312, India}

\author{M. C. Sharma}
\email{sharmamc2@gmail.com}
\address{Professor Department of Mathematics, N. R. E. C. College Khurja, Uttar Pradesh 203131, India}

\author{Amit Ujlayan}
\email{amitujlayan@gbu.ac.in}
\address{Department of Applied Mathematics, Gautam Buddha University, Greater Noida, Uttar Pradesh 201312, India}


\date{\today}
	
\begin{abstract}
This paper introduces and explores functions defined on \( H^* \)-normal spaces through the framework of \( H^* \)-open sets. We extend the concept of \( H^* \)-normality and investigate its connections with \( g \)-normal and classical normal spaces. Additionally, we generalize \( H^* \)-closed and \( H^* \)-generalized closed functions, analyzing their fundamental properties. Several characterizations of \( H^* \)-normal spaces are established, along with preservation theorems that highlight the structural significance of these functions.
\end{abstract}
\maketitle
\section*{Introduction} 
The concept of normality holds a central position in topology, both for its foundational significance and for its utility in addressing problems across theoretical and applied domains. A major focus of ongoing research has been the characterization of normality via combinations of weaker separation and covering properties. Over the years, numerous generalizations of open and closed sets have been proposed to support such investigations. For instance, Stone \cite{stone1937applications} introduced the concept of regular-open sets, which was later extended by Levine through the introduction of semi-open sets in 1963 \cite{levine1963semi}, and subsequently \( g \)-closed sets in 1970 \cite{levine1970generalized}.
In 2015, Nidhi Sharma, Poonam Sharma, and collaborators introduced the notion of \( H^* \)-normal spaces \cite{Tomar2004normal}, stimulating further research on generalized closed sets and separation axioms within topological structures. Continuing in this direction, Hamant Kumar and Tomar N. (2024) proposed the concept of \( Q^* \)-normal spaces, a new form of generalized normality designed to capture weaker separation conditions \cite{Neeraj2024normal}. In the same year, Neeraj Kumar Tomar, Fahed Zulfeqarr, and M. C. Sharma\cite{tomar2025scnormalspacesfunctions} introduced the concept of \( SC^* \)-normal spaces along with certain related functions, offering further insights into weak forms of normality and expanding the study of separation axioms in topological spaces.

This paper is structured into four main sections, each addressing a distinct aspect of the study:

\textbf{Section 1.} This section reviews several generalized set types, including \( H^* \)-closed, \( gH^* \)-open, \( H^*g \)-closed, and related classes. We provide illustrative examples and analyze the interrelations among these sets. In addition, we revisit the definition of \( H^* \)-normal spaces, offering examples and presenting a fundamental equivalence theorem.

\textbf{Section 2.} We introduce and examine the concepts of strongly \( H^* \)-open, strongly \( H^* \)-closed, and almost \( H^* \)-irresolute functions. Several results are established that link these functions to the structure of \( H^* \)-normal spaces.

\textbf{Section 3.} This section focuses on various generalizations of closed functions, including \( gH^* \)-closed, \( H^*g \)-closed, quasi \( H^* \)-closed, \( H^* \)-\( H^*g \)-closed, \( H^* \)-\( gH^* \)-closed, and almost \( gH^* \)-closed functions. We explore fundamental properties and provide supporting theorems.

\textbf{Section 4.} We conclude with preservation theorems and characterizations related to \( H^* \)-normal spaces, supported by theorems, lemmas, and corollaries that highlight their structural significance.

\section{Preliminaries and Notations}

Unless stated otherwise, all spaces discussed in this work are considered to be topological spaces without assuming any particular separation axioms. Let \( \phi : (X, \tau) \rightarrow (Y, \sigma) \) $\&$ \( \psi : (Y, \sigma) \rightarrow (Z, \rho) \) ( \( \phi : X \rightarrow Y \) $\&$ \( \psi : Y \rightarrow Z \)) be functions between topological spaces. In this context, \( \phi \) denotes a function from the topological space \( (X, \tau) \) to \( (Y, \sigma) \), $\&$ \( \psi \) denotes a function from \( (Y, \sigma) \) to \( (Z, \rho) \). For any subset \( A \subseteq X \), the \emph{closure} of \( A \) is denoted by \( \overline{A} \) or \( \operatorname{cl}(A) \), and the \emph{interior} of \( A \) is denoted by \( \mathring{A} \) or \( \operatorname{int}(A) \).

We now define some basic notions which will be used throughout. For a good understanding of them, readers are referred to see \cite{hatir1996decomposition,levine1963semi, Sundaram2000,njȧstad1965some,rodrigo2012international,stone1937applications}.

\subsection{Definition:}

A subset \( A \) of a topological space \( X \) is said to be:

\begin{enumerate}
    \item \textbf{regular closed} \cite{stone1937applications}(briefly \(r\)-closed) if \( A = \operatorname{cl}(\operatorname{int}(A)) \).
    \item \textbf{semi-closed} \cite{levine1963semi} (briefly \(s\)-closed) if \( \operatorname{int}(\operatorname{cl}(A)) \subseteq A \).
    \item \textbf{\( w \)-closed} \cite{Sundaram2000} if \( \operatorname{cl}(A) \subseteq U \), whenever \( A \subseteq U \) and \( U \) is semi-open in \( X \).
    \item \textbf{\( \alpha \)-closed} \cite{njȧstad1965some} if \( \operatorname{cl}(\operatorname{int}(\operatorname{cl}(A))) \subseteq A \).
    \item \textbf{\( \alpha^* \)-set} \cite{hatir1996decomposition} if \( \operatorname{int}(\operatorname{cl}(\operatorname{int}(A))) = \operatorname{int}(A) \).
    \item \textbf{\( C \)-set} \cite{hatir1996decomposition} if \( A = U \cap V \), where \( U \) is an open set and \( V \) is an \( \alpha^* \)-set in \( X \).
    \item \textbf{\( h \)-closed} \cite{rodrigo2012international} if \( s\text{-}\operatorname{cl}(A) \subseteq U \), whenever \( A \subseteq U \) and \( U \) is \( w \)-open in \( X \).
    \item \textbf{\( gh \)-closed} \cite{rodrigo2012international} if \( h\text{-}\operatorname{cl}(A) \subseteq U \), whenever \( A \subseteq U \) and \( U \) is \( h \)-open in \( X \).
    \item \textbf{regular-\( h \)-open} \cite{rodrigo2012international} if there exists a regular open set \( U \) such that \( U \subseteq A \subseteq h\text{-}\operatorname{cl}(U) \).
    \item \textbf{\( rgh \)-closed} \cite{rodrigo2012international} if \( h\text{-}\operatorname{cl}(A) \subseteq U \), whenever \( A \subseteq U \) and \( U \) is regularly \( h \)-open in \( X \).
    \item \textbf{\( hCg \)-closed} \cite{rodrigo2012international} if \( h\text{-}\operatorname{cl}(A) \subseteq U \), whenever \( A \subseteq U \) and \( U \) is a \( C \)-set in \( X \).
\end{enumerate}

\noindent
The complement of a \(r\)-closed set (or $s$, \( w \), \( \alpha \), \( h \), \( gh \), \( rgh \), and \( hCg \)) closed sets is called a \textbf{$r$-open} set (or \textbf{$s$}, \( w \), \( \alpha \), \( h \), \( gh \), \( rgh \), and \( hCg \)) \textbf{open sets}, respectively.

\noindent
The \textbf{\( h \)-closure} (or \textbf{\( s \)-closure}) of a subset \( A \) of a topological space, denoted by \( \overline{A}^h \) (or \( \overline{A}^s \)), is the smallest \( h \)-closed (respectively, \( s \)-closed) set that contains \( A \). Formally, it is defined as
\[
\overline{A}^h = \bigcap \{ F \subseteq X \mid A \subseteq F \text{ and } F \text{ is } h\text{-closed} \},
\]
and
\[
\overline{A}^s = \bigcap \{ F \subseteq X \mid A \subseteq F \text{ and } F \text{ is } s\text{-closed} \}.
\]

\noindent
The \textbf{\( h \)-interior} (or \textbf{\( s \)-interior}) of a set \( A \), denoted \( \mathring{A}^h \) (or \( \mathring{A}^s \)), refers to the greatest \( h \)-open (respectively, \( s \)-open) subset fully contained in \( A \). It is constructed as the union of all \( h \)-open (or \( s \)-open) sets that lie entirely within \( A \).

\subsection{Definition}

A subset \( A \) of a topological space \( X \) is said to be \textbf{\( H^* \)-closed} \cite{Tomar2004normal} if \( h\text{-}\operatorname{cl}(A) \subseteq U \) whenever \( A \subseteq U \) and \( U \) is \( hCg \)-open in \( X \). The complement of a \( H^* \)-closed set is called \textbf{\( H^* \)-open}.\\
\noindent
\noindent
The \textbf{\( H^* \)-closure} of a set \( A \subseteq X \), denoted by \( \overline{A}^{H^*} \) or \( H^*\text{-cl}(A) \), is the intersection of all \( H^* \)-closed sets containing \( A \). That is,
\[
H^*\text{-cl}(A) = \bigcap \{ C \subseteq X \mid A \subseteq C \text{ and } C \text{ is } H^*\text{-closed} \}.
\]

\noindent
The \textbf{\( H^* \)-interior} of \( A \), written as \( \mathring{A}^{H^*} \) or \( H^*\text{-int}(A) \), is the union of all \( H^* \)-open sets contained in \( A \), i.e.,
\[
H^*\text{-int}(A) = \bigcup \{ O \subseteq X \mid O \subseteq A \text{ and } O \text{ is } H^*\text{-open} \}.
\]

The family of all \( H^* \)-open sets (respectively, \( H^* \)-closed sets, $r$-open sets, $r$-closed sets, $s$-open sets, $\&$ $s$-closed sets) in a space \( X \) is denoted by \( H^*O(X) \) (respectively, \( H^*C(X) \), \( RO(X) \), \( RC(X) \), \( SO(X) \), $\&$ \( SC(X) \)).

\subsection{Definition}

Let \( A \) be a subset of a topological space \( (X, \tau) \). Then:

\begin{enumerate}
    \item \textbf{\( g \)-closed} \cite{levine1970generalized} if \( \operatorname{cl}(A) \subseteq U \) whenever \( A \subseteq U \) and \( U \in \tau \).
    
    \item \textbf{generalized \( H^* \)-closed} \cite{Tomar2004normal} (briefly, \( gH^* \)-closed) if \( H^*\text{-}\operatorname{cl}(A) \subseteq U \) whenever \( A \subseteq U \) and \( U \) is \( H^* \)-open.
    
    \item \textbf{\( H^* \)-generalized closed} (briefly, \( H^*g \)-closed) if \( H^*\text{-}\operatorname{cl}(A) \subseteq U \) whenever \( A \subseteq U \) and \( U \in X\).
    
    \item \textbf{regular-\( H^* \)-open} if there exists a regular open set \( U \) such that \( U \subseteq A \subseteq H^*\text{-}\operatorname{cl}(U) \).
\end{enumerate}

The complement of a \( g \)-closed (respectively, \( gH^* \)-closed, \( H^*g \)-closed) set is called \( g \)-open (respectively, \( gH^* \)-open, \( H^*g \)-open).  
The complement of a regular-\( H^* \)-open set is called regular-\( H^* \)-closed.

\subsection{Remark}

The following sequence of implications illustrates the hierarchy among various generalized closed sets in topological spaces:

\vspace{5mm}

\begin{center}
\begin{tabular}{c}
\text{closed} \(\longrightarrow\) \(\alpha\)-closed \(\longrightarrow\) \(h\)-closed \(\longrightarrow\) \(H^*\)-closed \(\longrightarrow\) \(gh\)-closed \(\longrightarrow\) \(rgh\)-closed \\
\hspace{7.5cm} \(\downarrow\) \hspace{2.cm} \(\downarrow\) \\
\hspace{7.2cm} \(gH^*\)-closed \(\longrightarrow\) \(rgH^*\)-closed
\end{tabular}
\end{center}

\vspace{5mm}

It is important to note that while each implication holds in the forward direction, the reverse implications are not generally valid. This can be demonstrated with suitable counterexamples.

\subsection{Example}

Let \( X = \{p, q, r, s\} \) be a topological space with the topology  
\[ 
\tau = \{\emptyset, \{p\}, \{q\}, \{p, q\}, \{p, q, r\}, X\}. 
\]  
Consider the subset \( A = \{r\} \subseteq X \).  

Then \( A \) is an \( h \)-closed set as well as an \( H^* \)-closed set in \( X \), but it is not closed in the topological sense, since its complement \( \{p, q, s\} \) is not open in \( \tau \).  

This example demonstrates that a set can be \( h \)-closed and \( H^* \)-closed without being closed in the usual topological sense.

\subsection{Example}

Let \( X = \{p, q, r, s, t\} \) be a topological space with the topology  
\[
\tau = \{\emptyset, \{p\}, \{s\}, \{t\}, \{p, s\}, \{p, t\}, \{s, t\}, \{p, s, t\}, X\}.
\]  
Consider the subset \( A = \{p, s, t\} \subseteq X \).  

The set \( A \) is both \( rgh \)-closed and \( rgH^* \)-closed in \( X \), but neither \( gh \)-closed nor \( gH^* \)-closed.  

Example illustrates, \( rgh \)-closedness $\&$ \( rgH^* \)-closedness do not imply \( gh \)-closedness or \( gH^* \)-closedness, emphasizing the distinctness of these generalized closure properties.

\subsection{Example}

Let \( X = \{p, q, r, s\} \) be a topological space with the topology  
\[
\tau = \{\emptyset, \{p\}, \{q\}, \{p, q\}, \{p, q, r\}, X\}.
\]  
Consider the subset \( A = \{r\} \subseteq X \).  

Thus \( A \) is a \( gH^* \)-closed set in \( X \), but not a closed set in the standard topological sense, since its complement \( \{p, q, s\} \) is not open in \( \tau \).  

This example demonstrates that a set can be \( gH^* \)-closed without being topologically closed.

 \subsection{Definition}
 A topological space $X$ is said to be normal\cite{tietze1923beitrage,zaitsev1968some} \textbf{(resp. $g$-normal\cite{malathi2017pre}, $H^*$-normal}) if for any pair of disjoint closed sets $A$ and $B$, there exist open (resp. $g$-open, $H^*$-open) sets $U$ and $V$ such that $A\subset U$ and $B\subset V$.

\subsection{Remark}

In a topological space \( X \), the following hierarchy of normality conditions holds:

\vspace{3.5mm}

\noindent
\textbf{normal} \(\Rightarrow\) \textbf{\( g \)-normal} \(\Rightarrow\) \textbf{\( H^* \)-normal}

\vspace{3.5mm}

\noindent
It is important to note that the converse of each implication does not generally hold, as counterexamples can be constructed to demonstrate their failure.

\subsection{Example}

Let the topological space \( X = \{p, q, r, s\} \) equipped with the topology
\[
\tau = \left\{ \emptyset, \{p\}, \{q\}, \{p, q\}, \{r, s\}, \{p, r, s\}, \{q, r, s\}, X \right\}.
\]

Consider the sets \( A = \{s\} \) $\&$ \( B = \emptyset \), which are disjt., closed. $\exists$ disjt., open sets \( U = \{p, r, s\} \) $\&$ \( V = \{q\} \), s.t. \( A \subseteq U \) $\&$ \( B \subseteq V \).  

Thus, \( X \) satisfies the conditions for being normal, \( g \)-normal, $\&$ \( H^* \)-normal, since every open set in \( X \) is also an \( H^* \)-open set.

\subsection{Theorem} 
For a topological space \( X \), the following statements are equivalent:

\begin{enumerate}
    \item \( X \) is \( H^* \)-normal.
    
    \item For every pair of open sets \( Q \) $\&$ \( R \) with \( Q \cup R = X \), there $\exists$ \( H^* \)-closed sets \( D \subseteq Q \) $\&$ \( E \subseteq R \) s.t. \( D \cup E = X \).
    
    \item For all closed set \( J \subseteq X \) $\&$ every open set \( K \) with \( J \subseteq K \), there $\exists$ an \( H^* \)-open set \( Q \) s.t.
    \[
    J \subseteq Q \subseteq H^*\text{-}\operatorname{cl}(Q) \subseteq K.
    \]
\end{enumerate}
\textbf{Proof.} 
\noindent
\textbf{(1) \( \Rightarrow \) (2)}: Let \( Q \) $\&$ \( R \) be open sets s.t. \( X = Q \cup R \). The sets \( X \setminus Q \) $\&$ \( X \setminus R \) are disjt., closed sets. Since \( X \) is \( H^* \)-normal, there \( \exists \) disjt., \( H^* \)-open sets \( Q_1 \) $\&$ \( R_1 \) s.t. \( X \setminus Q \subseteq Q_1 \) $\&$ \( X \setminus R \subseteq R_1 \). Define \( D = X \setminus Q_1 \) $\&$ \( E = X \setminus R_1 \). Then, \( D \) $\&$ \( E \) are \( H^* \)-closed sets satisfying \( D \subseteq Q \), \( E \subseteq R \), $\&$ \( D \cup E = X \).

\vspace{1.7mm}

\noindent
\textbf{(2) \( \Rightarrow \) (3)}: Let \( J \) be a closed $\&$ \(J\subseteq K\) $\&$ \( K \) an open set. The sets \( X \setminus J \), \( K \) are open $\&$ union is \(J\cup K = X \). By \textbf{(2)}, \( \exists \) \( H^* \)-closed sets \( M_1 \) $\&$ \( M_2 \) s.t. \( M_1 \subseteq X \setminus J \) $\&$ \( M_2 \subseteq K \), with \( M_1 \cup M_2 = X \). Thus, \( J \subseteq X \setminus M_1 \) $\&$ \( X \setminus K \subseteq X \setminus M_2 \), implying that \( (X \setminus M_1) \cap (X \setminus M_2) = \emptyset \). Let \( Q = X \setminus M_1 \) $\&$ \( R = X \setminus M_2 \). Then, \( Q \) $\&$ \( R \) are disjt., \( H^* \)-open sets s.t. \( J \subseteq Q \subseteq X \setminus R \subseteq K \). Since \( X \setminus R \) is \( H^* \)-closed, we have \( H^* \text{-}\operatorname{cl}(Q) \subseteq X \setminus R \), $\&$ hence \( J \subseteq Q \subseteq H^* \text{-}\operatorname{cl}(Q) \subseteq K \).

\vspace{1.7mm}

\noindent
\textbf{(3) \( \Rightarrow \) (1)}: Let \( J_1 \) $\&$ \( J_2 \) be two disjt., closed sets. Let \( K = X \setminus J_2 \), so \( J_2 \cap K = \emptyset \). Since \( J_1 \subseteq K \), $\&$ \( K \) is open, by \textbf{(3)} there \( \exists \) an \( H^* \)-open set \( Q \) s.t. \( J_1 \subseteq Q \subseteq H^* \text{-}\operatorname{cl}(Q) \subseteq K \). It follows that \( J_2 \subseteq X \setminus H^* \text{-}\operatorname{cl}(Q) = R \), where \( R \) is \( H^* \)-open. Hence, \( Q \cap R = \emptyset \), showing that \( J_1 \) $\&$ \( J_2 \) are separated by disjt., \( H^* \)-open sets \( Q \) $\&$ \( R \).

Therefore, \( X \) is \( H^* \)-normal.

\section{Functions Associated with \( H^* \)-Normal Spaces}

\subsection{Definition} A function \( \phi: X \to Y \) is called:

\begin{enumerate}
    \item \textbf{R-map} \cite{carnahan1974} if for every regular open set \( V \) in \( Y \), the preimage \( \phi^{-1}(V) \) is \(r\)-open in \( X \).
    \item \textbf{completely continuous} \cite{arya1974strongly} if for every open set \( V \) in \( Y \), the preimage \( \phi^{-1}(V) \) is regular open in \( X \).
    \item \textbf{rc-continuous} \cite{jankovic1985note} if for every regular closed set \( F \) in \( Y \), the preimage \( \phi^{-1}(F) \) is regular closed in \( X \).
    \item \textbf{strongly \( H^* \)-open} if for each \( U \in H^*O(X) \), it holds that \( \phi(U) \in H^*O(Y) \).
    \item \textbf{strongly \( H^* \)-closed} if for each \( U \in H^*C(X) \), it holds that \( \phi(U) \in H^*C(Y) \).
    \item \textbf{almost \( H^* \)-irresolute} if for every \( x\) in \(X\) $\&$ every \( H^* \)-nbd \( V \) of \( \phi(x) \), the closure\\ \( H^* \)-closure of the preimage \( \phi^{-1}(V) \), i.e., \( H^* \)-cl(\( \phi^{-1}(V) \)), is a \( H^* \)-nbd of \( x \).
\end{enumerate}

\subsection{Theorem}
A mapping \(\phi: X\) to \(Y\) is strong., \( H^* \)-closed iff, for every subset \( D \subseteq Y \) $\&$ every \( H^* \)-open set \( Q \subseteq X \) containing \( \phi^{-1}(D) \), $\exists$ a \( H^* \)-open set \( S \subseteq Y \) s.t. \( D \subseteq S \) $\&$ \( \phi^{-1}(S) \subseteq Q \).

\begin{flushleft}
\textbf{Proof.}
\textbf{(\( \Rightarrow \))} Suppose \( \phi \) is strongly \( H^* \)-closed. Let \( D \subseteq Y \), $\&$ let \( Q \in H^*O(X) \) be s.t. \( \phi^{-1}(D) \subseteq Q \). Define
\[
S = Y \setminus \phi(X \setminus Q).
\]
Since \( X \setminus Q \in H^*C(X) \), $\&$ \( \phi \) maps \( H^* \)-closed sets to \( H^* \)-closed sets, we have \( \phi(X \setminus Q) \in H^*C(Y) \). Hence, \( S \in H^*O(Y) \). Clearly, \( D \subseteq S \) $\&$ \( \phi^{-1}(S) \subseteq Q \), as required.
\end{flushleft}

\begin{flushleft}
\textbf{(\( \Leftarrow \))} Conversely, assume the given condition holds. Let \( K \in H^*C(X) \), so \( X \setminus K \in H^*O(X) \). Define \( D = Y \setminus \phi(K) \), then
\[
\phi^{-1}(D) \subseteq X \setminus K.
\]
By assumption, $\exists$ a \( H^* \)-open set \( S \subseteq Y \) s.t. \( D \subseteq S \) $\&$ \( \phi^{-1}(S) \subseteq X \setminus K \). Thus, \( \phi(K) \supseteq Y \setminus S \), which implies \( \phi(K) \in H^*C(Y) \). Therefore, \( \phi \) is strong.,\( H^* \)-closed.
\end{flushleft}

\subsection{Lemma}
For a mapping \(\phi: X\) to \(Y\), Each of the following assertions is equivalent:

\begin{enumerate}
    \item \( \phi \) is almost-\( H^* \)-irresolute.
    \item For every \( S \in H^*O(Y) \), we have
    \[
    \phi^{-1}(S) \subseteq H^*\text{-}int\big(H^*\text{-}cl(\phi^{-1}(S))\big).
    \]
\end{enumerate}

\subsection{Theorem}
Let \(\phi: X\) to \(Y\) be a mapping between topological spaces. Then \( \phi \) is almost \( H^* \)-irresolute iff
\[
H^*\text{-}\operatorname{cl}(Q) \subseteq H^*\text{-}\operatorname{cl}(\phi(Q)) \quad \text{for all } Q \in H^*\mathcal{O}(X).
\]

\textbf{Proof.}
\noindent
(\( \Rightarrow \)) Suppose \( \phi \) is almost \( H^* \)-irresolute, let \( Q \in H^*\mathcal{O}(X) \). Assume that \( y \notin H^*\text{-}\operatorname{cl}(\phi(Q)) \). Then there $\exists$ an \( H^* \)-open set \( S \subseteq Y \) s.t. \( y \in S \) $\&$ \( S \cap \phi(Q) = \emptyset \). 

This implies \( \phi^{-1}(S) \cap Q = \emptyset \). Since \( \phi \) is almost \( H^* \)-irresolute, it follows that:
\[
\phi^{-1}(S) \subseteq H^*\text{-}\operatorname{int}(H^*\text{-}\operatorname{cl}(\phi^{-1}(S))).
\]

Hence,
\[
H^*\text{-}\operatorname{int}(H^*\text{-}\operatorname{cl}(\phi^{-1}(S))) \cap H^*\text{-}\operatorname{cl}(Q) = \emptyset,
\]
which further implies
\[
\phi^{-1}(S) \cap H^*\text{-}\operatorname{cl}(Q) = \emptyset \quad \Rightarrow \quad S \cap \phi(H^*\text{-}\operatorname{cl}(Q)) = \emptyset.
\]

Therefore, \( y \notin \phi(H^*\text{-}\operatorname{cl}(Q)) \), since \( y \) was arbitrary, we conclude:
\[
H^*\text{-}\operatorname{cl}(Q) \subseteq H^*\text{-}\operatorname{cl}(\phi(Q)).
\]

\medskip

\noindent
(\( \Leftarrow \)) Now suppose that for every \( Q \in H^*\mathcal{O}(X) \), the inclusion \( H^*\text{-}\operatorname{cl}(Q) \subseteq H^*\text{-}\operatorname{cl}(\phi(Q)) \) holds. Let \( S \in H^*\mathcal{O}(Y) \), and define \( M = X \setminus H^*\text{-}\operatorname{cl}(\phi^{-1}(S)) \), which belongs to \( H^*\mathcal{O}(X) \). 

By the hypothesis, we have:
\[
\phi(H^*\text{-}\operatorname{cl}(M)) \subseteq H^*\text{-}\operatorname{cl}(\phi(M)).
\]

Now observe:
\[
\begin{aligned}
X \setminus H^*\text{-}\operatorname{int}(H^*\text{-}\operatorname{cl}(\phi^{-1}(S))) 
&= H^*\text{-}\operatorname{cl}(M) \\
&\subseteq \phi^{-1}(H^*\text{-}\operatorname{cl}(\phi(M))) \\
&\subseteq \phi^{-1}(H^*\text{-}\operatorname{cl}(Y \setminus S)) = X \setminus \phi^{-1}(S).
\end{aligned}
\]

Thus, \( \phi^{-1}(S) \subseteq H^*\text{-}\operatorname{int}(H^*\text{-}\operatorname{cl}(\phi^{-1}(S))) \), and so by the characterization of almost \( H^* \)-irresolute functions, \( \phi \) is almost \( H^* \)-irresolute.

\subsection{Theorem.}
Let \(\phi: X\) to \(Y\) be a surjective mapping that is cont., strongly \(H^*\)-open, $\&$ almost \(H^*\)-irresolute. If \(X\) is an \(H^*\)-normal space, then \(Y\) is also \(H^*\)-normal.

\textbf{Proof.}
Let \( D \subseteq Y \) be a closed subset $\&$ \( E \subseteq Y \) an open subset s.t. \( D \subseteq E \). Since \( \phi \) is cont., the preimage \( \phi^{-1}(D) \) is closed in \( X \), $\&$ \( \phi^{-1}(E) \) is open in \( X \), with \( \phi^{-1}(D) \subseteq \phi^{-1}(E) \).

Given that \( X \) is \( H^* \)-normal, there $\exists$ a \( H^* \)-open set \( Q \subseteq X \) satisfying:
\[
\phi^{-1}(D) \subseteq Q \subseteq H^*\text{-cl}(Q) \subseteq \phi^{-1}(E).
\]

Applying \( \phi \) to the above chain of inclusions, we obtain:
\[
D \subseteq \phi(Q) \subseteq \phi(H^*\text{-cl}(Q)) \subseteq E.
\]

Since \( \phi \) is strongly \( H^* \)-open, it follows that \( \phi(Q) \in H^*O(Y) \). Also, because \( \phi \) is almost \( H^* \)-irresolute, we have:
\[
H^*\text{-cl}(\phi(Q)) \subseteq E.
\]

Hence, there $\exists$ a \( H^* \)-open set \( \phi(Q) \) in \( Y \) s.t.:
\[
D \subseteq \phi(Q) \subseteq H^*\text{-cl}(\phi(Q)) \subseteq E.
\]

By the definition of \( H^* \)-normality, this confirms that \( Y \) is \( H^* \)-normal.

\subsection{Theorem.}
Let \(\phi: X\) to \(Y\) be a cont., and strongly \(H^*\)-closed surjection, where \(X\) is an \(H^*\)-normal space. Then the space \(Y\) is also \(H^*\)-normal.

\begin{flushleft}
\textbf{Proof.} Suppose \( M_1 \) $\&$ \( M_2 \) be two disjt., closed subsets of \( Y \). Then their inverse images \( \phi^{-1}(M_1) \) $\&$ \( \phi^{-1}(M_2) \) are closed in \( X \). Since \( X \) is \( H^* \)-normal, there $\exists$ disjt., \( H^* \)-open sets \( Q \) $\&$ \( S \) in \( X \) s.t.:
\[
\phi^{-1}(M_1) \subseteq Q \quad \text{and} \quad \phi^{-1}(M_2) \subseteq S.
\]

Because \( \phi \) is strongly \( H^* \)-closed, there $\exists$ \( H^* \)-open sets \( D \) $\&$ \( E \) in \( Y \) s.t.:
\[
M_1 \subseteq D, \quad M_2 \subseteq E, \quad \phi^{-1}(D) \subseteq Q, \quad \text{and} \quad \phi^{-1}(E) \subseteq S.
\]

Since \( Q \) $\&$ \( S \) are disjt., it follows that \( \phi^{-1}(D) \cap \phi^{-1}(E) = \emptyset \), which implies \( D \cap E = \emptyset \). Therefore, \( D \) $\&$ \( E \) are disjt., \( H^* \)-open sets containing \( M_1 \) $\&$ \( M_2 \), respectively. Hence, \( Y \) is \( H^* \)-normal.
\end{flushleft}

\subsection{Theorem.}
Let a mapping \(\phi: X\) to \(Y\) be a surjective function. Then, \( \phi \) is almost \( H^*g \)-closed iff for every subset \( T \subseteq Y \) $\&$ for each $r$-open set \( Q \subseteq X \) s.t. \( \phi^{-1}(T) \subseteq Q \), there $\exists$ a \( H^*g \)-open set \( R \subseteq Y \) s.t. \( T \subseteq R \) $\&$ \( \phi^{-1}(R) \subseteq Q \).

\vspace{1.7mm}

\textbf{Proof.}
\textbf{(Necessity)}: Assume \( \phi \) is almost \( H^*g \)-closed. Let \( T \subseteq Y \) be a subset, $\&$ let \( Q \subseteq X \) be a $r$-open set s.t. \( \phi^{-1}(T) \subseteq Q \). Define the set \( R = Y \setminus \phi(X \setminus Q) \). By construction, \( R \) is \( H^*g \)-open in \( Y \), $\&$ \( R \) contains \( T \). Furthermore, since \( \phi^{-1}(R) \subseteq Q \), the condition is satisfied.

\vspace{1.5mm}

\textbf{(Sufficiency)}: Suppose the condition holds. Let \( D \subseteq X \) be a $r$-closed set. Then, \( X \setminus D \) is a $r$-open set. Define \( E = Y \setminus \phi(D) \). Clearly,
\[
\phi^{-1}(E) \subseteq X \setminus D.
\]
Since \( X \setminus D \) is $r$-open, by assumption, there $\exists$ a \( H^*g \)-open set \( R \subseteq Y \) s.t.:
\[
E \subseteq R \quad \text{and} \quad \phi^{-1}(R) \subseteq X \setminus D.
\]
Thus, we obtain:
\[
\phi(D) \supseteq Y \setminus R \quad \Rightarrow \quad \phi(D) = Y \setminus R.
\]
This shows that \( \phi(D) \) is \( H^*g \)-closed in \( Y \). Hence, \( \phi \) is almost \( H^*g \)-closed.

\section{Extensions of \( H^* \)-Closed Function Concepts}

\subsection{Definition.} A function \(\phi: X \rightarrow Y\) is said to be:

\begin{enumerate}
    \item \textbf{\( H^* \)-closed}\cite{Tomar2004normal} if \( \phi(D) \) is \( H^* \)-closed in \( Y \) for every closed set \( D \subseteq X \).
    \item \textbf{\( H^*g \)-closed}\cite{Tomar2004normal} if \( \phi(D) \) is \( H^*g \)-closed in \( Y \) for every closed set \( D \subseteq X \).
    \item \textbf{\( gH^* \)-closed} if \( \phi(D) \) is \( gH^* \)-closed in \( Y \) for every closed set \( D \subseteq X \).
    \item \textbf{quasi \( H^* \)-closed} if \( \phi(D) \) is closed in \( Y \) for every \( D \in H^*C(X) \).
    \item \textbf{\( H^* \)-\( H^*g \)-closed} if \( \phi(D) \) is \( H^*g \)-closed in \( Y \) for every \( D \in H^*C(X) \).
    \item \textbf{\( H^* \)-\( gH^* \)-closed} if \( \phi(D) \) is \( gH^* \)-closed in \( Y \) for every \( D \in H^*C(X) \).
    \item \textbf{almost \( gH^* \)-closed}\cite{Tomar2004normal} if \( \phi(D) \) is \( gH^* \)-closed in \( Y \) for every \( D \in RC(X) \), where \( RC(X) \) denotes the regular closed sets of \( X \).
  \item \textbf{$H^*$-$gH^*$-continuous} if $\phi^{-1}(K)$ is $gH^*$-closed in $X$ for every $K\in H^*C(Y)$.
  \item \textbf{$H^*$-irresolute} if $\phi^{-1}(V)\in H^*O(X)$ for every $V\in H^*O(Y)$.
\end{enumerate}

\subsection{Theorem.}
Let \(\phi: X\) to \(Y\) $\&$ \( \psi: Y\) to \(Z\) be functions. Then:

\begin{itemize}
    \item[(1)] The mapping \( \psi \circ \phi: X \) to \(Z\) is \( H^* \)-\( gH^* \)-closed if \( \phi \) is \( H^* \)-\( gH^* \)-closed $\&$ \( \psi \) is cont., and \( H^* \)-\( gH^* \)-closed.
    
    \item[(2)] The mapping \( \psi \circ \phi: X \) to \(Z\) is \( H^* \)-\( gH^* \)-closed if \( \phi \) is strongly \( H^* \)-closed $\&$ \( \psi \) is \( H^* \)-\( gH^* \)-closed.
    
    \item[(3)] The mapping \( \psi \circ \phi: X \) to \(Z\) is \( H^* \)-\( gH^* \)-closed if \( \phi \) is quasi \( H^* \)-closed $\&$ \( \psi \) is \( gH^* \)-closed.
\end{itemize}

\subsection{Theorem.}
Let \( \phi: X\) to \(Y\) $\&$ \( \psi: Y \) to \(Z\) be mappings s.t. the mapping \( \psi \circ \phi: X \) to \(Z \) is \( H^* \)-\( gH^* \)-closed. If \( \phi \) is an \( H^* \)-irresolute surjection, then \( \psi \) is also \( H^* \)-\( gH^* \)-closed.

\begin{flushleft}
\textbf{Proof.}
Consider any set \( K \in H^*C(Y) \), i.e., a \( H^* \)-closed subset of \( Y \). Since \( \phi \) is an \( H^* \)-irresolute surjective mapping, its preimage \( \phi^{-1}(K) \) belongs to \( H^*C(X) \). Given that the composition \( \psi \circ \phi \) is \( H^* \)-\( gH^* \)-closed, we have:
\[
(\psi \circ \phi)(\phi^{-1}(K)) = \psi(K),
\]
which implies \( \psi(K) \) is \( gH^* \)-closed in \( Z \). Hence, \( \psi \) maps \( H^* \)-closed sets in \( Y \) to \( gH^* \)-closed sets in \( Z \), proving that \( \psi \) is \( H^* \)-\( gH^* \)-closed.
\end{flushleft}

\subsection{Lemma.} 
A function \(\phi: X\) to \(Y\) is called \( H^* \)-\( gH^* \)-closed iff for every subset \( E \subseteq Y \) and every \( H^* \)-open set \( Q \subseteq X \) s.t. \( \phi^{-1}(E) \subseteq Q \), there $\exists$ a \( gH^* \)-open set \( S \subseteq Y \) s.t. \( E \subseteq S \) $\&$ \( \phi^{-1}(S) \subseteq Q \).

\begin{flushleft}
\textbf{Proof.} 
\textbf{(\( \Rightarrow \))} Assume that \( \phi \) is \( H^* \)-\( gH^* \)-closed. Let \( E \subseteq Y \) be any subset, let \( Q \in H^*O(X) \) be s.t. \( \phi^{-1}(E) \subseteq Q \). Define \( S = Y \setminus \phi(X \setminus Q) \). Then \( S \) is a \( gH^* \)-open set in \( Y \), and since \( E \subseteq S \) $\&$ \( \phi^{-1}(S) \subseteq Q \), the condition is satisfied.

\textbf{(\( \Leftarrow \))} Conversely, suppose that for each \( E \subseteq Y \) and each \( H^* \)-open set \( Q \subseteq X \) with \( \phi^{-1}(E) \subseteq Q \), there $\exists$ a \( gH^* \)-open set \( S \subseteq Y \) s.t. \( E \subseteq S \) $\&$ \( \phi^{-1}(S) \subseteq Q \). Let \( W \in H^*C(X) \) be a closed set in \( X \). Consider the set \( Y \setminus \phi(W) \). Then \( \phi^{-1}(Y \setminus \phi(W)) \subseteq X \setminus W \), which is \( H^* \)-open. By the assumption, there $\exists$ a \( gH^* \)-open set \( S \subseteq Y \) s.t. \( Y \setminus \phi(W) \subseteq S \) $\&$ \( \phi^{-1}(S) \subseteq X \setminus W \). This implies \( \phi(W) \supseteq Y \setminus S \), so \( \phi(W) = Y \setminus S \), which shows that \( \phi(W) \) is \( gH^* \)-closed in \( Y \). Therefore, \( \phi \) is \( H^* \)-\( gH^* \)-closed.
\end{flushleft}

\subsection{Theorem.} 
Let \(\phi: X\) to \(Y\) be a cont., $\&$ \( H^* \)-\( gH^* \)-closed mapping. Then, for every \( gH^* \)-closed subset \( N \subseteq X \), the image \( \phi(N) \) is \( gH^* \)-closed in \( Y \).

\begin{flushleft}
\textbf{Proof.} 
Let \( N \subseteq X \) be a \( gH^* \)-closed set, $\&$ let \( S \subseteq Y \) be an open set s.t. \( \phi(N) \subseteq S \). Since \( \phi \) is cont., \( \phi^{-1}(S) \) is open in \( X \) and contains \( N \). Thus, \( H^*\text{-}\mathrm{cl}(N) \subseteq \phi^{-1}(S) \), which implies:
\[
\phi(H^*\text{-}\mathrm{cl}(N)) \subseteq S.
\]
Because \( \phi \) is \( H^* \)-\( gH^* \)-closed $\&$ \( H^*\text{-}\mathrm{cl}(N) \in H^*C(X) \), it follows that:
\[
gH^*\text{-}\mathrm{cl}(\phi(N)) \subseteq gH^*\text{-}\mathrm{cl}(\phi(H^*\text{-}\mathrm{cl}(N))) \subseteq S.
\]
Hence, \( \phi(N) \) is \( gH^* \)-closed in \( Y \).
\end{flushleft}

\subsection{Remark.} 
It holds that every \( H^* \)-irresolute function is \( H^* \)-\( gH^* \)-cont., However, the converse does not necessarily hold.

\subsection{Theorem.} 
A mapping \(\phi: X\) to \(Y\) is said to be \( H^* \)-\( gH^* \)-cont., iff the preimage \( \phi^{-1}(S) \) is \( gH^* \)-open in \( X \) for every \( S \in H^*O(Y) \).

\subsection{Theorem.} 
Let \(\phi: X\) to \(Y\) be an \( H^* \)-\( gH^* \)-cont., mapping. Then for every \( gH^* \)-closed subset \( W \subseteq Y \), the preimage \( \phi^{-1}(W) \) is \( gH^* \)-closed in \( X \).

\begin{flushleft}
\textbf{Proof.} 
Let \( W \subseteq Y \) be a \( gH^* \)-closed set, and let \( Q \in H^*O(X) \) be an \( H^* \)-open set s.t. \( \phi^{-1}(W) \subseteq Q \). Define the set \( S = Y \setminus \phi(X \setminus Q) \). Then \( S \) is a \( gH^* \)-open subset of \( Y \), with \( W \subseteq S \) $\&$ \( \phi^{-1}(S) \subseteq Q \).

Since \( W \subseteq S \), it follows that \( H^*\text{-}\mathrm{cl}(W) \subseteq S \). Taking preimages, we get:
\[
\phi^{-1}(W) \subseteq \phi^{-1}(H^*\text{-}\mathrm{cl}(W)) \subseteq \phi^{-1}(S) \subseteq Q.
\]

As \( \phi \) is \( H^* \)-\( gH^* \)-cont., the set \( \phi^{-1}(H^*\text{-}\mathrm{cl}(W)) \) is \( gH^* \)-closed in \( X \). Hence,
\[
gH^*\text{-}\mathrm{cl}(\phi^{-1}(W)) \subseteq gH^*\text{-}\mathrm{cl}(\phi^{-1}(H^*\text{-}\mathrm{cl}(W))) \subseteq Q.
\]

Therefore, \( \phi^{-1}(W) \) is \( gH^* \)-closed in \( X \), completing the proof.
\end{flushleft}

\subsection{Corollary.}
Let \(\phi: X\) to \(Y\) be a closed $\&$ \( H^* \)-irresolute mapping. Then for every \( gH^* \)-closed set \( W \subseteq Y \), the preimage \( \phi^{-1}(W) \) is \( gH^* \)-closed in \( X \).

\subsection{Theorem.}
Let \(\phi: X\) to \(Y\) be an open, bijective, and \( H^* \)-\( gH^* \)-cont., mapping. Then for every \( gH^* \)-closed set \( W \subseteq Y \), the preimage \( \phi^{-1}(W) \) is \( gH^* \)-closed in \( X \).

\begin{flushleft}
\textbf{Proof.} 
Let \( W \subseteq Y \) be \( gH^* \)-closed, let \( Q \subseteq X \) be an open set s.t. \( \phi^{-1}(W) \subseteq Q \). Since \( \phi \) is open and surjective, we have:
\[
W = \phi(\phi^{-1}(W)) \subseteq \phi(Q),
\]
$\&$ \( \phi(Q) \) is open in \( Y \). Hence, \( H^*\text{-cl}(W) \subseteq \phi(Q) \).

Because \( \phi \) is injective, it follows that:
\[
\phi^{-1}(W) \subseteq \phi^{-1}(H^*\text{-cl}(W)) \subseteq \phi^{-1}(\phi(Q)) = Q.
\]

Since \( \phi \) is \( H^* \)-\( gH^* \)-cont., the preimage \( \phi^{-1}(H^*\text{-cl}(W)) \) is \( gH^* \)-closed in \( X \). Therefore,
\[
gH^*\text{-cl}(\phi^{-1}(W)) \subseteq gH^*\text{-cl}(\phi^{-1}(H^*\text{-cl}(W))) \subseteq Q.
\]
Thus \( \varphi^{-1}(W) \) is \( gH^* \)-closed in \( X \), as required.
\end{flushleft}

\subsection{Theorem.}
Let \( \phi: X\) to \(Y\) $\&$ \( \psi: Y\) to \(Z\) be mappings s.t. the mapping \( \psi \circ \phi: X \) to \(Z\) is \( H^* \)-\( gH^* \)-closed. If \( \psi \) is an open, bijective, and \( H^* \)-\( gH^* \)-cont., mapping, then \( \phi \) is \( H^* \)-\( gH^* \)-closed.

\begin{flushleft}
\textbf{Proof.}
Let \( W \in H^*C(X) \). Since \( \psi \circ \phi \) is \( H^* \)-\( gH^* \)-closed, it follows that
\[
(\psi \circ \phi)(W) = \psi(\phi(W))
\]
is \( gH^* \)-closed in \( Z \). Because \( \psi \) is a bijection, we have
\[
\phi(W) = \psi^{-1}(\psi(\phi(W))).
\]

Moreover, since \( \psi \) is \( H^* \)-\( gH^* \)-cont., $\&$ open, the inverse image of a \( gH^* \)-closed set under \( \psi^{-1} \) is also \( gH^* \)-closed in \( Y \) by \textbf{Theorem 3.10}. Hence, \( \phi(W) \in gH^*C(Y) \), $\&$ so \( \phi \) is \( H^* \)-\( gH^* \)-closed.
\end{flushleft}

\subsection{Theorem.}
Let \(\phi: X\) to \(Y\) $\&$ \( \psi: Y\) to \(Z\) be mappings s.t. the mapping \( \psi \circ \phi: X\) to \(Z\) is \( H^* \)-\( gH^* \)-closed. If \( \psi \) is a closed injection $\&$ \( H^* \)-\( gH^* \)-cont., then \( \varphi \) is \( H^* \)-\( gH^* \)-closed.

\begin{flushleft}
\textbf{Proof.}
Suppose \( W \in H^*C(X) \). Since \( \psi \circ \phi \) is \( H^* \)-\( gH^* \)-closed, the image \( (\psi \circ \phi)(W) = \psi(\phi(W)) \) is \( gH^* \)-closed in \( Z \). As \( \psi \) is injective, it follows that
\[
\phi(W) = \psi^{-1}(\psi(\phi(W))).
\]

Now, since \( \psi \) is \( H^* \)-\( gH^* \)-cont., and closed, by \textbf{Theorem 3.11,} the inverse image of a \( gH^* \)-closed set under \( \psi \) is \( gH^* \)-closed in \( Y \). Therefore, \( \phi(W) \in gH^*C(Y) \), which implies that \( \phi \) is \( H^* \)-\( gH^* \)-closed.
\end{flushleft}

\section{Preservation Theorems and Characterizations of \( H^* \)-Normal Spaces}

\subsection{Theorem.} 
Let \(\phi: X\) to \(Y\) be a cont., quasi \( H^* \)-closed surjection, $\&$ suppose that \( X \) is \( H^* \)-normal. Then, \( Y \) is normal.

\begin{flushleft}
\textbf{Proof.} 
Suppose \( G_1 \) $\&$ \( G_2 \) be two disjt., closed subsets of \( Y \). If \( \phi \) is cont., the preimages \( \phi^{-1}(G_1) \) $\&$ \( \phi^{-1}(G_2) \) are disjt., closed subsets of \( X \). Since \( X \) is \( H^* \)-normal, there $\exists$ disjt., \( H^* \)-open sets \( Q_1 \) $\&$ \( Q_2 \) in \( X \) s.t. \( \phi^{-1}(G_1) \subseteq Q_1 \) $\&$ \( \phi^{-1}(G_2) \subseteq Q_2 \).

Next, define the sets \( S_1 = Y \setminus \phi(X \setminus Q_1) \) $\&$ \( S_2 = Y \setminus \phi(X \setminus Q_2) \). These sets are open in \( Y \), and we have the following conditions:
\[
G_1 \subseteq S_1, \quad G_2 \subseteq S_2, \quad \phi^{-1}(S_1) \subseteq Q_1, \quad \phi^{-1}(S_2) \subseteq Q_2.
\]
Since \( Q_1 \cap Q_2 = \emptyset \) $\&$ \( \phi \) is surjective, it follows that \( S_1 \cap S_2 = \emptyset \).

Thus, \( Y \) is normal.
\end{flushleft}

\subsection{Lemma.}
A subset \( D \) of a topological space \( X \) is \( gH^* \)-open iff for every closed set \( F \) s.t. \( F \subseteq A \), it holds that \( F \subseteq H^*\)-int(\( A \)).

\subsection{Theorem.} 
Let \(\phi: X\) to \(Y\) be a closed \( H^* \)-\( gH^* \)-cont., injection, $\&$ suppose that \( Y \) is \( H^* \)-normal. Then \( X \) is \( H^* \)-normal.

\begin{flushleft}
\textbf{Proof.}  
Suppose closed set \( G_1 \) $\&$ \( G_2 \) be disjt.,\( X \). If \( \phi \) is a closed injection, the images \( \phi(G_1) \) $\&$ \( \phi(G_2) \) are disjt., closed sets in \( Y \). By the \( H^* \)-normality of \( Y \), there $\exists$ disjt., open sets \( S_1 \) $\&$ \( S_2 \) in \( H^*O(Y) \) s.t. \( \phi(G_i) \subseteq S_i \) for \( i = 1, 2 \).

Since \( \phi \) is \( H^* \)-\( gH^* \)-cont., the preimages \( \phi^{-1}(S_1) \) $\&$ \( \phi^{-1}(S_2) \) are disjt., \( gH^* \)-open sets in \( X \), $\&$ \( G_i \subseteq \phi^{-1}(S_i) \) for \( i = 1, 2 \).

Next, define \( Q_i = H^*\text{-}int(\phi^{-1}(S_i)) \) for \( i = 1, 2 \). These sets \( Q_i \) are elements of \( H^*O(X) \), \( G_i \subseteq Q_i \), and \( Q_1 \cap Q_2 = \emptyset \).

Thus, \( X \) is \( H^* \)-normal.
\end{flushleft}

\subsection{Corollary.}
Let \(\phi: X\) to \(Y\) be a closed \( H^* \)-irre., injection. If \( Y \) is \( H^* \)-normal, then \( X \) is \( H^* \)-normal.

\begin{flushleft}
    \textbf{Proof.} This follows from the fact that an \( H^* \)-irresolute function is \( H^* \)-\( gH^* \)-continuous. Since normality is preserved under closed \( H^* \)-\( gH^* \)-continuous injections, the result follows directly.
\end{flushleft}

\subsection{Lemma.}
A function \(\phi: X\) to \(Y\) is almost \( gH^* \)-closed iff for every subset \( E \subseteq Y \) $\&$ every $r$-open set \( Q \in RO(X) \) s.t. \( \phi^{-1}(E) \subseteq Q \), there $\exists$ a \( gH^* \)-open set \( S \subseteq Y \) s.t. \( E \subseteq S \) $\&$ \( \phi^{-1}(S) \subseteq Q \).

\subsection{Lemma.} 
Let \(\phi: X\) to \(Y\) be an almost \( gH^* \)-closed function. Then for every closed set \( G \subseteq Y \) and for every $r$-open set \( Q \in RO(X) \) s.t. \( \phi^{-1}(G) \subseteq Q \), there $\exists$ a set \( S \in H^*O(Y) \) s.t. \( G \subseteq S \) $\&$ \( \phi^{-1}(S) \subseteq Q \).

\section{Conclusion}

This paper presents the introduction and analysis of a new class of spaces called \( H^* \)-normal spaces, which are defined using \( H^* \)-open sets. We investigate the connections between traditional normality and \( H^* \)-normality, demonstrating that \( H^* \)-normality is a more generalized or weaker form of normality. Additionally, we provide various characterizations of \( H^* \)-normal spaces and discuss the properties of \( gH^* \) and \( H^*g \)-closed functions within this framework. The findings in this study contribute valuable insights and methods that can assist researchers in further exploring normal spaces of this type. Moreover, the concepts and techniques developed here could be extended to other types of topological spaces, such as ordered topological spaces, bi-topologically ordered spaces, and fuzzy topological spaces, offering potential for broad applications in topological research.


\bibliographystyle{amsplain}
\bibliography{references}

\end{document}